\documentclass[12pt]{article}
\usepackage{amsfonts, amsmath, amssymb,latexsym,amsthm}
\usepackage[margin=1.3in]{geometry}
\usepackage[all]{xy}

\newtheorem{defn0}{Definition}[section]
\newtheorem{prop0}[defn0]{Proposition}
\newtheorem{thm0}[defn0]{Theorem}
\newtheorem{lemma0}[defn0]{Lemma}
\newtheorem{corollary0}[defn0]{Corollary}
\newtheorem{example0}[defn0]{Example}
\newtheorem{remark0}[defn0]{Remark}
\newtheorem{conjecture0}[defn0]{Conjecture}

\newenvironment{definition}{ \begin{defn0}}{\end{defn0}}
\newenvironment{proposition}{\bigskip \begin{prop0}}{\end{prop0}}
\newenvironment{theorem}{\bigskip \begin{thm0}}{\end{thm0}}

\newenvironment{corollary}{\bigskip \begin{corollary0}}{\end{corollary0}}
\newenvironment{example}{ \begin{example0}\rm}{\end{example0}}
\newenvironment{remark}{ \begin{remark0}\rm}{\end{remark0}}

\newcommand{\propref}[1]{Proposition~\ref{#1}}
\newcommand{\thmref}[1]{Theorem~\ref{#1}}

\newcommand{\corref}[1]{Corollary~\ref{#1}}
\newcommand{\exref}[1]{Example~\ref{#1}}

\def\max{{\bf m}}                   
\def\res{{\bf k}}                   
\def\rees{{\mathcal R}}             
\def\length{{\mathrm{ Length}}}
\def\depth{{\mathrm{ depth}}}
\def\grCM{{\mathrm{ gCM}}}
\def\agrCM{{\mathrm{ agCM}}}
\def\spec{{\bf Spec}}
\def\proj{{\bf Proj}}
\def\hir{{\bf Hir}}
\def\loc{{\bf Loc}}

\def\highfull#1#2#3{{\bf R}^{#1}#2_*{#3}}
\def\highd#1#2#3{{\bf R}^{#1}#2_*{\mathcal O}_{#3}}
\def\high#1#2#3{\length_R({\bf R}^{#1}#2_*{\mathcal O}_{#3})}
\def\HiRI+#1#2{H^{#2}_{\rees(#1)_{+}}\rees(#1)}

\def\oc{{\mathcal O_C}}

\title{  \bf \huge On the last Hilbert-Samuel coefficient  of isolated
singularities
}
\author{\large   Juan Elias
\thanks{Partially supported by  MTM2007-67493, MTM2010-20279-C02-01}
 }

\date{November 1, 2010}

\begin{document}

\maketitle

\section{Introduction}

In   \cite{Lip78} Lipman presented  a proof  of the existence of a
desingularization for any excellent surface.
The strategy of Lipman's proof  is based on the finiteness of the number $H(R)$,
defined as the supreme
of the second Hilbert-Samuel coefficient $e_2(I)$, where $I$ range the set of
normal $\max$-primary ideals of
 a Noetherian complete local ring  $(R,\max)$.
See    \cite[Theorem* of pag. 158, and Remark B of pag. 160]{Lip78}.
On the other hand Huckaba and Huneke proved  that if $I$ is a $\max$-primary  ideal
of a $d-$dimensional, $d\ge 2$, Cohen-Macaulay local ring $(R,\max)$ such that $I^n$ is integrally closed for $n \gg 0$, in particular if $I$ is normal,  then the associated graded ring of $R$ with respect to $I^n$ has  depth at least two for $n\gg 0$, \cite[Theorem 3.11]{HH99}.

The problem studied  in this paper is the extension of the result of Lipman on $H(R)$ to $\max-$primary
ideals $I$ of a $d$-dimensional Cohen-Macaulay ring $R$ such that $gr_{I^n}(R)$,  associated graded ring of $R$ with respect to $I^n$, is Cohen-Macaulay for $n\gg 0$.
We denote by  $\grCM(R)$ this set of ideals of $R$.

 We prove the following theorem:

\bigskip
\noindent
\textbf{Theorem 3.2}
\emph{
Let $(R,\max)$ be a $d-$dimensional Cohen-Macaulay local ring of dimension $d\ge 1$ essentially of finite type over a characteristic zero field $\res=R/\max$.
Assume that the closed point  $x \in X=\spec(R)$ is an isolated singularity.
For all $K\in   \grCM(R)$ and for all Hironaka  ideal $I\in \hir(R)$ it holds
$$
  0\le e_d(K)\le e_d(I)= p_g(X).
$$
In particular, if   $x\in X$ is a rational singularity then
$ e_d(K)=0$ for all  $K\in \grCM(R)$.
}

\bigskip
\noindent
 An Hironaka ideal of $R$  is an ideal $I$ such that the blow-up of $R$ centered at $I$
 is a desingularization of $X=\spec(R)$.
  We know from the main result of \cite[Main Theorem I, pag. 132]{Hir64a} that
  Hironaka ideals exist.
 Recall that all Hironaka ideal belongs to $\grCM(R)$, this is a consequence of  the version  of the  Grauert-Riemenschneider vanishing
 theorem due to Sancho de Salas, \propref{basicprop}.

 An example of Huckaba and Huneke shows that  cannot be considered in \thmref{ineq} weaker conditions on the depth of the associated graded ring with respect $I^n$,  and that cannot be  extended the result of Lipman to higher dimensions, see \exref{hhex}.
 In \propref{extremal} we relate the ideals such that their  last Hilbert-Samuel coefficient equals the geometric genus with the rational singularities.
In the last section we study the one-dimension case; in particular we explicitly construct an Hironaka ideal, \exref{notfix}.
We also relate the upper bound of $e_1(I)$ of the main theorem of this paper with
the previously obtained in \cite{Eli08} and \cite{RV05}.

\bigskip
\noindent {\sc Acknowledgments:}
We thank O. Villamayor for the useful comments and suggestions about the resolution of singularities.

\bigskip
\noindent {\sc Notations:}
Let $I$ be an ideal of $R$ we denote by
$\rees(I)=\oplus_{n\ge0} I^n t^n$  the Rees algebra of $I$, and  we denote by $gr_I(R)=\oplus_{n\ge
0}I^n/I^{n+1}$ the associated graded ring of $R$ with respect to
$I$.
If $I$ is a $\max-$primary ideal of $R$ we denote by
$h_{I}(n)=length_R(R/I^{n+1})$ the Hilbert-Samuel function of $I$.
Hence there exist integers $e_j(I) \in \mathbb Z$ such that
$$
p_{I}(X)=
e_0(I) \binom{X+d}{d}-  e_1(I)\binom{X+d-1}{d-1} + \cdots +  (-1)^d e_d(I)
$$
is the  Hilbert-Samuel polynomial of $I$, i.e. $h_{I}(n)=p_{I}(n)$
for $n \gg 0$. The integer $e_j(I)$  is the $j-$th normalized
Hilbert-Samuel coefficient of $I$, $j=0,\cdots , d$.
We set $h_R=h_{\max}$ and  $p_R=p_{\max}$.
We denote by $P_I$ the Poincar\'{e} series of an $\max$-primary ideal $I$,
i.e. the power series defined by $h_I$
$$
P_I(Z)=\sum_{n\ge 0} h_I(n)\; Z^n. $$
We know that $P_I$ is a rational function, with $a_i\in \mathbb Z$ and  $a_s\neq 0$,
$$
P_I=\frac{a_0+a_1 Z+\cdots +a_s Z^s}{(1-Z)^d}.
$$


\section{On the resolution of singularities}

Let $(R,\max)$ be a $d-$dimensional reduced Cohen-Macaulay local ring essentially of finite type over a characteristic zero field $\res=R/\max$.
We denote by $x$ the closed point of $X=\spec(R)$.
Hironaka  proved  that there exists an ideal $I\subset R$ such that,
\cite[Main Theorem I, pag. 132]{Hir64a},
\begin{enumerate}
\item[(i)] $V(I)=Sing(X)$,
\item[(ii)] the natural projection morphism $\pi  : \widetilde{X}=\proj(\rees (I))\longrightarrow X=\spec(R)$
is a resolution of sin\-gu\-larities of $X$, i.e. $\tilde{X}$ is non-singular and $\pi$ induces  a $\res$-scheme isomorphism
$$
\pi  : \widetilde{X}\setminus \pi^{-1}(Sing(X)) \longrightarrow X\setminus Sing(X).
$$
\end{enumerate}

\begin{definition}
An ideal $I$ of $R$ is an  Hironaka
ideal if $I$ satisfies the two above conditions.
We denote by $\hir(R)$ the set of Hironaka ideals of $R$.
\end{definition}

See \exref{notfix} for an explicit computation of an Hironaka ideal.

\medskip
In the next result we collect three basic properties of Hironaka ideals.

\begin{proposition}
\label{basicprop}
Let $I$ be an Hironaka ideal, then
\begin{enumerate}
\item[(i)] $I^n\in \hir(R)$, for  $n\ge 1$,
\item[(ii)] $gr_{I^n}(R)$ is Cohen-Macaulay for all $n\gg 0$,
\item[(iii)] if the closed point $x$ of $X$ is an isolated singularity then
 Hironaka's ideals are $\max$-primary.
\end{enumerate}
\end{proposition}
\begin{proof}
$(i)$ The result follows form the facts $V(I^n)=V(I)=Sing(X)$ and
that there exists a natural $\res$-scheme  isomorphism
$$
\proj(\rees (I))\cong \proj(\rees (I^n))
$$
induced by the degree zero graded morphism $\rees(I^n)\subset \rees(I)$, $n\ge 1$.

\noindent
$(ii)$ Follows from the version of the  Grauert-Riemenschneider vanishing theorem due to Sancho de Salas, \cite{San87b}.
See \cite[Theorem 4.3]{Lip94} for an extension
to all Cohen-Macaulay rings.

\noindent
$(iii)$ If the closed point $x$ of $X$ is an isolated singularity then $V(I)=Sing(X)=\{x\}$,
i.e. $I$ is an $\max$-primary ideal.
\end{proof}

\begin{remark}
\label{charpfails}
Recall that if the field $\res$ is of positive characteristic it is an open problem
to prove that  $\hir(R)$ is non empty.
On the other hand, the fact that Hironaka Ideal $I$ has a Cohen-Macaulay associated graded ring $gr_{I^n}$ for $n\gg0$ is
the key point of the results of this paper.
This a consequence of the version of the Grauert-Riemenschneider vanishing theorem due to Sancho de Salas, \cite{San87b}, as we quoted in the proof of \propref{basicprop}.
Recall that  Grauert-Riemenschneider vanishing theorem is true in characteristic zero case but not in the positive characteristic case as  Raynaud showed in \cite{Ray78}.
\end{remark}

\begin{definition}
If the closed point of $X$ is an  isolated singularity,
the geometric genus  of $X$ is
$p_g(X)=  \high{d-1}{\pi}{{Z}}$
for a (all)  singularity resolution
$\pi: Z\longrightarrow X$.
We say that $R$  is a rational singularity if $R$ is normal and
for all $i>0$ $\high{i}{\pi}{{Z}}=0$.
\end{definition}

It is well know that the above definitions do not depend on the resolution $Z$.
In fact from the Grothendieck spectral sequence of the composition  of functors and the vanishing result of the higher direct  images due to Hironaka we get
that $\high{d-1}{\pi}{{Z}}$ is independent of the resolution of singularities $\pi:Z \longrightarrow X$, see \cite[Proposition 2.1]{HO74}.
In the one-dimensional  case,
$\highd{0}{\pi}{{Z}}\cong \overline{R}/R$ where $\overline{R}$ is the integral closure of $R$
on its full ring of fractions.
In the literature the length of   $\overline{R}/R$ is also known as the singularity order
of $X=\spec(R)$.
In the one-dimensional case, rational means non-singular.

\begin{remark}
 Recall that if $R$ is Cohen-Macaulay of dimension $d\ge 2$ and the closed point of $X=\spec(R)$ is an  isolated singularity then $R$ is a reduced normal ring.
\end{remark}

\medskip
We assume that  $(R,\max)$ is  a $d-$dimensional Noetherian Cohen-Macaulay local ring of dimension $d\ge 1$.

\begin{definition}
We denote by $\grCM(R)$ (resp. $\agrCM(R)$)     the set of  $\max$-primary ideals $I$ of $R$
such that $gr_{I^n}(R)$ is Cohen-Macaulay (resp. $\depth(gr_{I^n}(R)) \ge d-1$) for $n\gg 0$.
\end{definition}

Notice that $\hir(R)\subset \grCM(R) \subset \agrCM(R)$, these inclusions are in general strict.
See \cite{Eli04} and its reference list for a detailed study of the ideals belonging to $\grCM(R)$ or $\agrCM(R)$.
Moreover, if $I$ is a $\max$-primary ideal it is known that
$$
\depth(gr_{I^n}(R))\ge 1
$$
for $n\gg 0$.
In particular, if $d\le 2$ then $\agrCM(R)$  agrees with the whole set of $\max$-primary ideals.
If  $d=1$ then $\agrCM(R)=\grCM(R)$, and this set agrees with the set of $\max$-primary ideals.

The relationship between the  Hilbert-Samuel coefficients  of $I^n$ and $I$   is the following.

\begin{proposition}
Let $(R,\max)$ be a Noetherian local ring and let $I$ be a $\max$-primary ideal.
Then it hold:
\begin{enumerate}
\item[(i)]  $e_d(I^n)=e_d(I)$, $n\ge 1$,
\item[(ii)] there exist  polynomials $f_i^j$  on $n$ of degree $d-j$,
$i=0,\cdots,d-1$, $j=0,\cdots, i$, such that
$$
e_i(I^n)=\sum_{j=0}^{i} f_i^j(n) e_j(I).
$$
\end{enumerate}
\end{proposition}

\medskip
\noindent
Since $e_d(I^n)= e_d(I)$ for all $n\ge 1$, in order to study the last Hilbert-Samuel
coefficient
of ideals of $\grCM(R)$, we may
assume that $gr_I(R)$ is Cohen-Macaulay.
In the following result we present some inequalities among the Hilbert-Samuel coefficient of $\max$-primary ideals
$I$ such that $gr_I(R)$ is Cohen-Macaulay.
 The same result holds for Hilbert filtrations
$\mathcal F$ such that its associated graded ring $gr_{\mathcal F}(R)$ is Cohen-Macaulay.
In \cite[Corollary 2]{Mar89} Marley gave several restrictions for the Hilbert-Samuel coefficients of ideals $I\in \agrCM(R)$.

\begin{proposition}
\label{hhc}
Let $(R,\max)$ be a $d-$dimensional Cohen-Macaulay local ring of dimension $d\ge 0$.
Let $I$ be an $\max$-primary ideal such that $gr_I(R)$ is Cohen-Macaulay with Poincar\'{e} series, $a_s\neq 0$,
$$
P_I=\frac{a_0+a_1 Z+\cdots +a_s Z^s}{(1-Z)^d}.
$$
Then the following conditions hold
\begin{enumerate}
\item[(i)] $s\le e_0(I) +d +1 -\length(I/I^2)\le e_0(I)$,
\item[(ii)] $e_i(I)=0$, for  $i=s+1,\cdots, d$,
\item[(iii)] $0\le (i+1) e_{i+1}(I) \le (s-i) e_i(I)$, for $i=0,\cdots, s$,
\item[(iv)] $0\le e_i(I)\le \binom{s}{i} e_0(I)$, for $i=0,\cdots, d$.
\end{enumerate}
\end{proposition}
\begin{proof}
We know that, \cite{EV91},
$$
e_i(I)=\sum_{j\ge i} \binom{j}{i} a_j
$$
for $i=0,\cdots, d$, from this we deduce  $(ii)$.
Since $gr_I(R)$ is Cohen-Macaulay we get that $a_i>0$ for all $i=0,\cdots, s$, and
$a_1=\length(I/I^2)-d$.
Then we have $e_i(I)\ge 0$ for all $i=0,\cdots,d$, and
$$
e_0(I)\ge a_0+ a_1+ (a_2+\cdots + a_s)\ge a_0 + a_1 + s-1,
$$
since $a_0=\length(R/I)$, we  obtain $(i)$.
An easy computation shows that
$$
(s-i) \binom{j}{i} \ge (i+1) \binom{j}{i+1}
$$
for $j=i,\cdots, s$. From this we get $(iii)$.
The inequalities of $(iv)$ follow from the previous ones.
\end{proof}

\medskip
\begin{example}
Let us assume that  $s=5$. From the last result we get, $e_i=e_i(I)$,
$$
 10 e_0 \ge 2 e_1 \ge  e_2 \ge  e_3 \ge 2 e_4 \ge  10 e_5\ge 0.
$$
Compare these inequalities with \cite[Corollary 3.5]{PUVV10}.
\end{example}

\bigskip
\section{Hilbert-Samuel coefficients of isolated singularities}

Let $(R,\max)$ be a $d-$dimensional Cohen-Macaulay local ring of dimension $d\ge 2$ essentially of finite type over a characteristic zero field $\res=R/\max$.
In this section we always assume that the closed point of $X=\spec(R)$ is an isolated singularity.

The first part of the next result it is well known, we include it for the sake of completeness.

\begin{proposition}
\label{upperbound}
$(i)$ Let $I$ be a $\max$-primary ideal of $R$, and let
$\pi:\widetilde{X}=\proj(\rees (I)) \longrightarrow X=\spec(R)$ be the projection, then
$$
\high{d-1}{\pi}{{\widetilde{X}}}=
\length_R(H^{d-1}(\widetilde{X}, {\cal O}_{\widetilde{X}}))=
\length_R(\HiRI+{I}{d}_0)< \infty.
$$

\noindent
$(ii)$
For all $I\in \agrCM(R)$ and for all $n\ge 1$ it holds
$$
e_d(I^n)=\high{d-1}{\pi}{{\widetilde{X}}}.
$$

\noindent
$(iii)$
For all Hironaka ideal $I\in \hir(R)$ and $n\ge 1$ it holds
$$
    e_d(I^n)= p_g(X).
$$
 In particular,  $e_d(I^n)$ is independent of the ideal $I\in \hir(R)$ and the integer $n\ge 1$.
\end{proposition}
\begin{proof}
$(i)$
Since $X$ is an affine scheme we get, \cite[Chap. III, Proposition 8.5]{Har97},
$$
 \highd{d-1}{\pi}{{\widetilde{X}}} \cong
 H^{d-1}(\widetilde{X}, {\cal O}_{\widetilde{X}})\widetilde\quad
$$
and by  \cite{Gro67} we have
$$
H^{d-1}(\widetilde{X}, {\cal O}_{\widetilde{X}}) \cong
\HiRI+{I}{d}_0.
$$
Since $x$ is an isolated singularity and $I$ is an $\max$-primary we get that
$\highd{d-1}{\pi}{{\widetilde{X}}}$ is supported in  $x$, so the lengths of $(i)$  are all finite and agree.

\noindent
$(ii)$
It is well known that  for all $n\ge 1$ it holds $e_d(I^n)=e_d(I)$
and there exists a  natural $\res$-scheme isomorphism
$$
\widetilde{X}=\proj(\rees (I))\cong \proj(\rees (I^n)).
$$
Hence we may assume that $\depth(gr_I(R))\ge d-1$.
From the Grothendieck-Serre formula, see
\cite[Proposition 3.1]{JV95}, we get
$$
\begin{array}{ccccc}
  e_d(I)&=& (-1)^d (p_I(-1)-h_I(-1))&= & \chi(\length_R(\HiRI+{I}{*}_0)) \\ \\
   & & &= & \length_R(\HiRI+{I}{d}_0).
\end{array}
 $$
From  $(i)$ we deduce $(ii)$.

\medskip
\noindent
$(iii)$ It is a consequence of $(ii)$ and \propref{basicprop}.
\end{proof}

\medskip
In the following result
we relate the last Hilbert-Samuel coefficient  $e_d(I)$ of a Hironaka ideal $I$ with the last Hilbert-Samuel coefficient $e_d(K)$ of a  $\max$-primary ideal $K\in \grCM(R)$.
In particular we extend \cite[(B) pag.160]{Lip78} to ideals $I$ such that $gr_{I^n}(R)$ is Cohen-Macaulay for $n\gg 0$.

\begin{theorem}
\label{ineq}
Let $(R,\max)$ be a $d-$dimensional Cohen-Macaulay local ring of dimension $d\ge 1$ essentially of finite type over a characteristic zero field $\res=R/\max$.
Assume that the closed point  $x \in X=\spec(R)$ is an isolated singularity.
For all $K\in   \grCM(R)$ and for all Hironaka  ideal $I\in \hir(R)$ it holds
$$
  0\le e_d(K)\le e_d(I)= p_g(X).
$$
In particular, if   $x\in X$ is a rational singularity then
$ e_d(K)=0$ for all  $K\in \grCM(R)$.
\end{theorem}
\begin{proof}
The one-dimensional case is studied in the last section, so we may assume that $d\ge 2$.
We set $\sigma: Z=\proj(\rees (K))\longrightarrow X$ the blowing-up of $X$ centered in $K$.
Since $e_d(K)=e_d(K^n)$ for all integer $n\ge 1$,
we may assume that $gr_K(R)$ is Cohen-Macaulay, \propref{basicprop} $(ii)$.

By the Principalization theorem and the existence of resolution of singularities we get that there exist an Hironaka ideal $I\subset R$ and
a commutative diagram, \cite[Section 1.9]{Kol07},
$$
\xymatrix{
\widetilde{X}=\proj(\rees (I))\ar[dr]^f \ar[dd]^{\pi}& \\
&Z=\proj(\rees (K)) \ar[dl]^\sigma\\
X=\spec(R)&\\
}
$$
Let  $E$ (resp. $\widetilde{E}$) be the exceptional divisor of
$\sigma$ (resp. $\pi$).
Since $gr_K(R)$ is Cohen-Macaulay,  from   \cite[Theorem 1.3]{San87b} we get
 that
 $$
 H^{i}_E(Z, {\cal O}_{Z})=0
 $$
 for $i<d$.
Hence we have the following exact sequence, \cite{Gro67},
$$
0=H^{d-1}_E(Z, {\cal O}_{Z})
\longrightarrow
H^{d-1}(Z, {\cal O}_{Z})
\stackrel{\rho_Z}{\longrightarrow}
H^{d-1}(Z-E, {\cal O}_{Z}),
$$
i.e. the restriction morphism $\rho_Z$ is  a monomorphism.
Since the morphisms $f$ and $\sigma$ induce  isomorphisms
$$
\widetilde{X} \setminus\widetilde{E}
\stackrel{f}{ \cong }
Z \setminus E
\stackrel{\sigma}{ \cong }
X\setminus \{x\}
,
$$
we have the commutative diagram
$$
\xymatrix{
H^{d-1}(Z, {\cal O}_{Z}) \ar[r]^{\rho_Z} \ar[d]^f & H^{d-1}(Z-E, {\cal O}_{Z}) \ar[d]^{\cong}\\
H^{d-1}(\widetilde{X}, {\cal O}_{\widetilde{X}}) \ar[r]^{\rho_{\widetilde{X}}}&
H^{d-1}(\widetilde{X}-\widetilde{E}, {\cal O}_{\widetilde{X}}).\\
}
$$
We know that $\rho_Z$ is a monomorphism, so  $f$ is also a
monomorphism.
By \propref{upperbound} we deduce
$$
e_d(K)=\length_R( H^{d-1}(Z, {\cal O}_{Z}))\le
\length_R( H^{d-1}(\widetilde{X}, {\cal O}_{\widetilde{X}}))=e_d(I)=p_g(X).
$$

\noindent
On the other hand, since $gr_K(R)$ is Cohen-Macaulay we have $e_d(K)\ge 0$,
\propref{hhc}.
\end{proof}

\medskip
\begin{remark}
A particular case of the second part of the last result is when $R$ is regular. In this case the result holds for a regular local ring $R$ without any restriction on  the field $\res$.
Without loss of generality we may assume that $k$ is infinite.
This is a well known result and a short proof could be  the following.
Let $(R,\max)$ be a $d$-dimensional, $d\ge 1$,  regular local ring with maximal ideal $\max$ and residue field $\res$.
Let $I$ be a $\max$-primary ideal such that $gr_{I^n}(R)$ is Cohen-Macaulay for some $n\ge 1$, we want to prove that $e_d(I)=0$.
If $d=1$ then $R$ is a DVR and the result is trivial.
Since $e_d(I)=e_d(I^n)$ for all integer $n\ge 1$,
we may assume that $gr_I(R)$ is Cohen-Macaulay.
Let $J$ be a minimal reduction of $I$.
Since $R$ is regular, from  Brian\c{c}on-Skoda theorem
we have  $I^{d-1}\subset J$, \cite{LT81} or  \cite[13.3.3]{SH-IC09}.
By \cite[Proposition 4.6]{HM97} we get $e_d(I)=0$.
\end{remark}

\begin{corollary}[{\cite[(B) pag.160]{Lip78}}]
\label{lipman}
Let $X=\spec(R)$ be a Cohen-Macaulay scheme of dimension $d=2$.
Assume that the closed point of $X$ is an isolated singularity.
Let $I$ be an $\max$-primary ideal.
 If $I^n$ is integrally closed for $n\gg 0$ then
$$
0 \le e_2(I) \le p_g(X).
$$
\end{corollary}
\begin{proof}
We know that $I\in \grCM(R)$,\cite{HH99}, so the claim follows from the last result.
\end{proof}

\begin{example}
\label{hhex}
The example of Huckaba and Huneke shows that \thmref{ineq} cannot be extended to ideals of $\agrCM(R)$ and that \corref{lipman} cannot be extended to higher dimensions, \cite[Theorem 3.11]{HH99}.
In fact, let us consider  $R=\res[X,Y,Z]_{(X,Y,Z)}$ with
$\res$ a characteristic zero field.
Notice that $R$ is a regular local ring of dimension $d=3$ defining a rational singularity $X=\spec(R)$.
Let $I$ the ideal of $R$ generated by $(X,Y,Z)^5$ and
$X^4, X(Y^3+Z^3), Y(Y^3+Z^3), Z(Y^3+Z^3)$.
Huckaba and Huneke proved that $I$ is a normal ideal, $\depth(gr_{I^n}(R))=2$
for all $n\ge 1$, i.e. $I\in \agrCM(R)\setminus \grCM(R)$,  and if $\widetilde{X}=\proj(\rees (I))$ then
$H^{2}(\widetilde{X}, {\cal O}_{\widetilde{X}})\neq 0$.
                         From \propref{upperbound} we get that
$$
p_g(X)=e_3(R)=0 < \length_R(H^{2}(\widetilde{X}, {\cal O}_{\widetilde{X}}))=e_3(I).
$$

Moreover, a standard computation gives us $e_3(I)$:
$$
p_{I}(X)=
76 \binom{X+2}{3}
- 48\binom{X+1}{2}
+ 4 \binom{X}{1}
-1
$$
i.e. $e_3(I)=1$.
\end{example}

\begin{remark}
Last results can not be extended to the low Hilbert-Samuel coefficients.
If $I\in \grCM(R)$   then $I^n\in \grCM(R)$ for all $n\ge 1$, but
$e_i(I^n)\rightarrow \infty$ when $n\rightarrow \infty$ for $i\neq d$,
see Section 2.
\end{remark}

\begin{remark}
\thmref{ineq} shows that $0\le e_d(I)\le p_g(X)$.
The extremal case $e_d(I)=0$ implies, under some conditions,
that the associated graded ring is Cohen-Macaulay.
Being $R$ a Cohen-Macaulay ring and  $I$ $\max$-primary we have $e_1(I)\ge e_0(I)-\length_R(R/I)\ge 0$, \cite{Nor60}, and $e_2(I)\ge 0$,  \cite{Sal81a}.
Huneke proved that if $e_1(I)=e_0(I)-\length_R(R/I)$ then $gr_{I^n}(R)$ is Cohen-Macaulay for all $n\ge 1$, \cite{Hun87}.
Narita proved that $e_2(I)=0$ if and only if $gr_{I^n}(R)$ is Cohen-Macaulay for all $n\gg 0$.
Marley gave an example of $\max$-primary ideal with $e_2(I)=0$ and that $gr_I(R)$ is not
Cohen-Macaulay.
On the other hand, Narita gave an example of an ideal of a Cohen-Macaulay ring
with $e_3(I)<0$, \cite{Nar63}.
See \cite[example 2]{Mar89} for  an  ideal of a regular local ring with a negative $e_3(I)$.
Itoh proved that if $I$ is normal then $e_3(I)\ge 0$, \cite{Ito92}.
Corso, Polini and Rossi proved that if $I^n$ is normal for some $n\gg 0$ then
$e_3(I)\ge 0$, and if $I^n$ is integrally closed for all $n\gg 0$ and $e_3(I)=0$ then
$gr_{I^n}(R)$ is Cohen-Macaulay for all $n\gg 0$,  \cite{CPR05}.
See \cite{Eli99} for the generalization of some of the above results to ideals
satisfying the second Vallabrega-Valla condition.
\end{remark}

\medskip
In the next result we look at the ideals with maximal $e_d(I)$.
See \cite[Proposition 2.7]{Mor84} for a related result with the
the first part of the next result.

\begin{proposition}
\label{extremal}
Let $K\in   \grCM(R)$ be an ideal such that $K^n$ is  integrally closed  for $n\gg 0$.
Assume that  $Z=\proj(\rees (K))$ has only isolated singularities.
Then $e_d(K)= p_g(X)$ if and only if
  $Z$ has only rational singularities.
\end{proposition}
\begin{proof}
Since $e_d(K)=e_d(K^n)$ for all integer $n\ge 1$ and $K\in   \grCM(R)$
we may assume that $gr_K(R)$ is Cohen-Macaulay, \propref{basicprop} $(ii)$.
Then  $Z=\proj(\rees (K))$ is a normal Cohen-Macaulay scheme.
Let $\widetilde{X}$ be  a desingularization of $Z$.
It is easy to prove that  $\widetilde{X}$  is also a desingularization of $X$.
Hence we have a commutative diagram
$$
\xymatrix{
\widetilde{X} \ar[dr]^f \ar[dd]^{\pi}& \\
&Z=\proj(\rees (K)) \ar[dl]^\sigma\\
X=\spec(R)&\\
}
$$

\noindent
Since $Z$ is a normal scheme we have
$\highd{0}{f}{\widetilde{X} }\cong {\mathcal O}_{Z}$, \cite{HO74}.
Being  $Z$  a Cohen-Macaulay scheme eventually with isolated singularities,
$\highd{i}{f}{{\widetilde{X} }}=0$, $i=1, \cdots, d-2$.
Hence from the Grothendieck spectral sequence we get, \cite[Theorem 11.3]{Rot79},
$$
0 \longrightarrow
\highfull{d-1}{\sigma}{{\mathcal O}_{Z}}
\longrightarrow
\highd{d-1}{\pi}{{\widetilde{X} }}
\longrightarrow
\highfull{0}{\sigma}{(\highd{d-1}{f}{\widetilde{X} })}
\longrightarrow
\highfull{d}{\sigma}{{\mathcal O}_{Z}}
\longrightarrow
0.
$$

\noindent
The arithmetical rank of $\rees(K)_{+}$ is $d=\dim(R)$, so
$$
\highfull{d}{\sigma}{{\mathcal O}_{Z}} \cong \HiRI+{K}{d+1}_0=0.
$$
From the last exact sequence and \propref{upperbound}  we get
$$
\begin{array}{ccl}
  \length_R(\highfull{0}{\sigma}{(\highd{d-1}{f}{\widetilde{X} })}) & =&
  \length_R(\highd{d-1}{\pi}{{\widetilde{X} }})-
  \length_R(\highfull{d-1}{\sigma}{{\mathcal O}_{Z}})\\ \\
      & = & p_g(X)-e_d(K).
\end{array}
$$

If     $Z$ has only rational singularities,
then we have $\highd{d-1}{f}{\widetilde{X} }=0$.
From the last equality we get  $e_d(K)= p_g(X)$.
If $e_d(K)= p_g(X)$  then
 $$
0= \highfull{0}{\sigma}{(\highd{d-1}{f}{\widetilde{X} })}\cong H^0(Z, \highd{d-1}{f}{\widetilde{X} })\widetilde \;\;  .
 $$
Since  the scheme $Z$ has only isolated singularities
we get  $\highd{d-1}{f}{\widetilde{X} }=0$, i.e.
 $Z$ has only rational singularities.
\end{proof}

See \cite{Mor84} and \cite{Mor87} for results related with the main result of this paper.

\section{The one-dimensional case}

Let $\widetilde{X}$ be the blowing-up of a finite union of reduced curve singularities $X=\spec(R)$ centered on one of its singular closed points, i.e. $R$ is a Cohen-Macaulay reduced semilocal ring.
Then   $\widetilde{X}$ has a finite number of non-singular closed points $P_1, \dots,
P_r$, and the germ $( \widetilde{X}, P_i)$ is a  curve singularity for
$i=1,\dots,r$.
It is known that the set $\{P_1, \dots, P_r \}$ is in
correspondence one-to-one with the set of closed  points of
${\mathbf {Proj}}(gr_{\max}(R))$.
We can iterate the process we get a sequence of blowing-ups

$$
\pi_i:X^{(i+1)} \longrightarrow X^{(i)}
$$

\noindent
$i\ge 0$, with $X^{(0)}=X$, $X^{(i+1)}$ is the blowing-up of $X^{i}$ centered at a singular point of
$X^{i}$, and $X^{(i)}$ is non-singular
for $i \ge r$.
The composition of all maps $\{\pi_i\}_{i=0,\cdots, r}$
is a resolution  $\pi:\widetilde{X}\longrightarrow X$.
It is well known that $\widetilde{X}\cong \spec(\overline{R})$,
where $\overline{R}$ is the integral closure of $R$ on its ring of fractions.
We denote by $\loc(X)$ the finite set of local rings ${\cal O}_{X^{(i)},q}$,
$q$ singular point of $X^{(i)}$,  appearing in the above resolution process. Recall that the family of Hilbert-Samuel polynomials
$p_{\oc}=e_o(\oc)(n+1)-e_1(\oc)$, with $\oc \in \loc(X)$ is univocally determined by $X$.

\begin{definition}
For a reduced curve singularity $X=\spec(R)$, the singularity
order of $X$ is the finite number
$$
\delta(X)=\length_R(\overline{R}/R),
$$
i.e. the geometric genus $p_g(X)$ of $X$.
\end{definition}

\noindent
From \cite{Nor59a} we have
$$
\delta(X)=\sum_{\oc \in \loc(X)} e_1(\oc).
$$
Recall that we can decompose $e_1(\oc)$ as sum of the micro-invariants of the ring extension $\oc \subset {\mathcal O_{Bl_{\max}(R)}}$, see \cite{Eli99}.

\begin{proposition}
\label{oned}
Let $X=\spec(R)$ be a one-dimensional reduced scheme.
Let $I$ be an $\max$-primary ideal, then
$$
0\le e_1(I)\le \delta(X)=p_g(X).
$$
If $e_1(I)=\delta(X)$ then $I$ is an Hironaka ideal.
\end{proposition}
\begin{proof}
We denote by $Bl_I(R)$ the semilocal  ring of the blow-up of $X=\spec(R)$ centered at the ideal $I$.
We know that the ring extensions $R\subset Bl_I(R)\subset \overline{R}$ are finite and that,  \cite{Lip71},
$$
0\le e_1(I)=\length_R(Bl_I(R)/R)\le \length_R(\overline{R}/R)=\delta(X).$$
 Hence we get the first part of the  result.
 If  $e_1(I)=\delta(X)$ then we get $Bl_I(R)=\overline{R} $, i.e. $I$  is an Hironaka ideal.
\end{proof}

\begin{example}
\label{notfix}
Let us consider the   one-dimensional Cohen-Macaulay local domain
$R_n=\res[x,y]_{(x,y)}/(y^2-x^n)$, $n\ge 8$;
we set $X_n=\spec(R_n)$.
The resolution process consists in $r=[n/2]$ blow-ups, and the singular points appearing in the process are all  of multiplicity two.
Since the $e_1$ of a double point is one we get $\delta(X)=r$.
Last result shows that for all $\max$-primary ideal $I\subset R_n$
$$
0\le e_1(I)\le \delta(X_n)=r.
$$
Let us consider the ideal $I_n$ of $R_n$ generated by $x^6$ and $x^2y$.
The Hilbert-Samuel polynomial of $I_n$ is
$p_{I_n}(t)= 12 t -4$, i.e.
$e_0(I_n)= 12$, and $e_1(I_n)=4$.
Hence, $I_n$ is an Hironaka ideal of $R_n$ if and only if
$n=8,9$.

The jacobian ideal $J_n=(y,x^{n-1})$ is not an Hironaka ideal because
$e_1(J_n)=1$.
Since $e_0(J_n)- \length(R_n/J_n)=1=e_1(J_n)$, we get that $gr_{J_n^t}(R_n)$ is Cohen-Macaulay for all $t\ge 1$, \cite{Hun87}.

On the other hand, in \cite{Eli08} and \cite{RV05} upper bounds
for $e_1$ are given.
In our case we get $e_1(I_n)=4 < \epsilon(I_n)=8$, bound of \cite{Eli08}, and $e_1(I_n)=4 < \rho(I_n)=59$, bound of \cite{RV05}.
See \cite[Proposition 2.2]{Eli08} for a comparison between $\epsilon$ and $\rho$.
\end{example}


\providecommand{\bysame}{\leavevmode\hbox to3em{\hrulefill}\thinspace}
\providecommand{\MR}{\relax\ifhmode\unskip\space\fi MR }
\providecommand{\MRhref}[2]{%
  \href{http://www.ams.org/mathscinet-getitem?mr=#1}{#2}
}
\providecommand{\href}[2]{#2}

\medskip
\noindent
Juan Elias\\
Departament d'\`Algebra i Geometria\\
Universitat de Barcelona\\
Gran Via 585, 08007 Barcelona, Spain\\
e-mail: {\tt elias@ub.edu}

\end{document}